\documentclass[numreferences]{kluwer}
\usepackage{amssymb,latexsym,eucal,mathrsfs}

\usepackage[curve,v2]{xypic}

\newcommand{\tensor}{\otimes}

\newcommand{\cork}{corank}

\newcommand{\picard}{Pic}

\newcommand{\wt}[1]{\widetilde{#1}}

\newcommand{\wtp}{\widetilde{\P^n}}
\newcommand{\wtv}{\widetilde{\varphi }}
\newcommand{\wts}{\widetilde{\Sigma}}

\newcommand{\blow}[2]{{\rm Bl}_{#2}{(#1})}

\newcommand{\ses}[3]{0\rightarrow#1\rightarrow#2
   \rightarrow#3\rightarrow0}

\newcommand{\TanX}{{\mathcal T}}
\newcommand{\SecX}{\Sigma}

\renewcommand{\H}{{\mathcal H}}
\newcommand{\F}{{\mathscr F}}
\newcommand{\D}{{\mathscr D}}
\newcommand{\E}{{\mathscr E}}

\newcommand{\T}{{\Theta}}
\renewcommand{\I}{{\mathscr I}}
\renewcommand{\L}{{\mathscr L}}
\renewcommand{\O}{{\mathcal O}}
\renewcommand{\P}{\dP}
\newcommand{\C}{\dC}

\newcommand{\Z}{\dZ}
\newcommand{\Q}{\dQ}

\newcommand{\G}{{\mathbb{G}}}
\newcommand{\of}{\, \mbox{\small{$\circ$}} \:}

\newenvironment{proof}{\par \medskip \noindent
{\sc Proof:}}{\nopagebreak \hfill $\Box$ \par \medskip}

\begingroup
\newtheorem{thm}{Theorem}[section]   
\newtheorem{cor}[thm]{Corollary}     
\newtheorem{lemma}[thm]{Lemma}         
\newtheorem{prop}[thm]{Proposition}  
\newtheorem{defn}[thm]{Definition}   
   
\endgroup

\newenvironment{rem}[2]{\refstepcounter{thm} \label{#2} 
\par \medskip \noindent {\em \bf #1 \thethm }}{\par \medskip}

\newcommand{\gmodh}[2]{\raisebox{.1cm}{#1}/\raisebox{-.1cm}{#2}} 

\begin{document}

\begin{article}
\begin{opening} 
\title{Results on Secant Varieties Leading to a Geometric 
Flip Construction}
\author{Peter \surname{Vermeire}\email{petvermi@math.okstate.edu}} 
\runningauthor{Peter Vermeire}
\institute{Oklahoma State University}
\begin{ao}
Department of Mathematics,
Oklahoma State University, 
Stillwater OK 74078
\end{ao}

\begin{abstract}
We study the relationship between the equations defining a
projective variety and properties of its secant varieties.
In particular, we use 
information about the syzygies among the defining equations to 
derive smoothness and normality statements about
$SecX$ and also to obtain information about linear systems on the 
blow up of projective space along 
a variety $X$.  We use these results to geometrically construct, for varieties 
of arbitrary dimension, a flip first described in the case of curves 
by M. Thaddeus via Geometric Invariant Theory.  
\end{abstract}
\keywords{syzygies, flips, vanishing theorems, secant varieties}
\classification{1991 Mathematics Subject Classification}{Primary
14E05; Secondary 13D02, 14F17, 14D20}

\end{opening}

\section{Introduction}

Let $X\subset \P^n$ be a projective variety, scheme theoretically
defined by homogeneous polynomials $F_0,\ldots, F_s$ of degree $d$.
The $F_i$ induce a rational map $\varphi:\P^n\dashrightarrow\P^s$
defined off $X$.  Equivalently, $\varphi$ is determined by a linear system
$V\subset \Gamma(\P^n,\O(d))$ with base scheme $X$. 
One resolves $\varphi$
to a morphism $\wtv:\wtp\rightarrow \P^s$ by blowing up $\P^n$ along $X$.  
$\varphi$ and $\wtv$ have been studied in a number of contexts, including
Cremona transformations \cite{ck1},\cite{ck2},\cite{es},\cite{hks}, linear 
systems on the blow up of projective space \cite{gl}, and in somewhat 
greater generality in the form of projections from subvarieties \cite{bhss}.

We will be most concerned with the case $d=2$.  It is
well known that if $\mathscr{L}$ is an ample line bundle on a variety
$X$, then $X$ is ideal theoretically defined by quadrics for all
embeddings induced by sufficiently large multiples of $\mathscr{L}$.
Explicit examples include smooth curves embedded by line bundles of
degree at least $2g+2$ and canonical curves with Clifford index 
at least $2$.

We begin with the observation that, in the case $d=2$, if $L\subset\P^n$
is a secant line to $X$ then $\varphi$ collapses $L$ to a point as the
restriction of the $F_i$ to $L$ forms a system of quadrics with a 
base scheme of length two, determining a unique quadric.  This raises
the natural question: When is $\varphi$ an embedding off $SecX$?
Corollary~\ref{weakembedding} shows this is the case if a condition slightly 
weaker than Green's condition $(N_2)$ is imposed.  In fact, in
Theorem~\ref{embed} we prove a more general statement (for 
arbitrary $d$) about when $\wtv$ is an embedding off the
proper transform of an appropriate secant variety.  We also give
in Corollary~\ref{littlevanishing} an application to the vanishing of
the cohomology of powers of ideal sheaves.

Section~\ref{secresults} is concerned with the structure of
$SecX$ and its proper transform in the case $d=2$.  
The main results in this 
section are Theorem~\ref{getabundle} and its Corollary~\ref{alliswell}
where it is shown that the restriction of $\wtv$ to the
proper transform of the secant variety
is a $\P^1$-bundle over the length two Hilbert scheme of $X$,
extending a result in \cite{bertram1} to varieties of
arbitrary dimension.  We also give (Remark~\ref{deficientsecants}) 
a geometric
criterion for determining the dimension of the secant
variety to a variety satisfying Green's condition $(N_2)$.

In Section~\ref{flipconst}, we obtain a partial
answer to a question raised in \cite[\S 1]{bertram4} by 
constructing a flip first described in the case of smooth
curves by M. Thaddeus \cite{thad} via GIT.
The construction proceeds in
several stages, with the end result summarized in 
Theorem~\ref{flip}.  The reader should note that these are not
$K_X$-flips in the sense of the Minimal Model Program.  It is shown in
\cite{bertram4} that the flips constructed by Thaddeus are, in
fact, log flips; however we do not address that question here.
Again (Corollary~\ref{secondvanishing}) we give a simple application to 
the cohomology of ideal sheaves.

Aside from gaining an understanding of the geometry of secant
varieties and how this geometry relates to syzygies, a practical
goal of the flip construction is the derivation of 
vanishing theorems for the groups
$H^i(\P^n,\I_X^a(b))$ (Cf. \cite{bel},\cite{wahl}).  
Specifically, in \cite{thad},
Thaddeus obtains vanishing theorems on the flipped
spaces via Kodaira vanishing as he is able to identify
the ample cone on each space.  In a related direction,
Bertram \cite{bertram3} uses a generalization of
Kodaira vanishing to prove vanishing theorems directly on the space
$\wtp$, deduced from the existence of log canonical divisors.
As discussed in \cite{bertram3}, a combination of the two
techniques should reveal the strongest results.  The construction of
further flips is taken up in \cite{verm-flip2}, and 
the question of vanishing theorems in \cite{verm-vanish}.

{\bf Notation:} We will decorate a projective variety $X$ as follows: 
$X^d$ is the $d^{th}$ cartesian product of $X$; $S^dX$ is
$Sym^dX=\gmodh{$X^d$}{$S_d$}$, the $d^{th}$ symmetric product of $X$;
and $\H^dX$ is $Hilb^d(X)$, the Hilbert Scheme of zero dimensional
subschemes of $X$ of length $d$.  Recall (Cf. \cite{gott}) that if $X$ is
a smooth projective variety then $\H^dX$ is also projective, and is smooth
if either $dim~X\leq 2$ or $d\leq 3$.

If $V$ is a $k$-vector space, $\P(V)$ is the space of
1-dimensional {\it quotients} of $V$.  We work
throughout over the field $k=\C$ of complex numbers.
We use the terms locally free sheaf (resp. invertible sheaf) and vector 
bundle (resp. line bundle) interchangeably.
If $D\subset X$ is a Cartier divisor, then the associated invertible
sheaf is denoted $\O_X(D)$.  We conform to the convention that 
products of line bundles corresponding to explicit divisors are written
additively, while other products are written multiplicatively, e.g. 
$(\L\tensor \O_X(D))^{\tensor n}\cong \L^{\tensor n}\tensor \O_X(nD)$.  
A line bundle $\L$ on $X$ is {\em nef} if $\L.C\geq 0$ for every 
irreducible curve $C\subset X$.  A line bundle $\L$ is {\em big} if 
$\L^{\tensor n}$ induces a birational map for all $n\gg 0$.
The term {\em conic} is used to mean a quadric hypersurface
in some projective space.

{\bf Acknowledgments:}
I would like to thank the referee for many useful comments and
suggestions.  I would also I like to thank the following people for
their helpful  
conversations and communications:
Aaron Bertram, Lawrence Ein, Anthony Geramita, 
Klaus Hulek, Shreedhar Inamdar, Vassil Kanev, S\'andor Kov\'acs, Mario Pucci, 
M. S. Ravi, Michael Schlessinger, and Michael Thaddeus.
I would lastly like to thank Jonathan Wahl for his 
guidance, enthusiasm and limitless patience in teaching me the 
subject of Algebraic Geometry and in supervising my dissertation,
from which most of these results are derived.

\section{Condition $(K_d)$}

We begin with the situation $d=2$ from the Introduction and
establish a useful technical result:
\begin{prop}\label{embeddingequiv}
Let $V\subseteq \Gamma(\P^n,\O(2))$ be a linear system with base scheme
$X$, with induced map $\varphi:\P^n\dashrightarrow\P(V)$. 
Then the following are equivalent:
\begin{enumerate}
\item $\varphi$ is an embedding off $SecX$
\item If $L\subset \P^n$ is a line not intersecting $X$, $L\not\subset
SecX$, then the natural restriction map $r_L:V\rightarrow \Gamma(L,\O_L(2))$
is surjective
\end{enumerate}
\end{prop}

\begin{proof}
Assume $\varphi$ is an embedding off $SecX$, and let $L\not\subset
SecX$ be a line not intersecting $X$.
Then the
restriction of $\varphi$ to $L$ is base point free, hence 
$\cork(r_L)\leq 1$.  If $\cork(r_L)=1$, then $\varphi|_L$ is a 
ramified double cover of $\P^1$, contradiction the assumption that
$\varphi$ is an embedding off $SecX$.

Conversely, choose a length two 
subscheme $Z\subset \P^n\setminus SecX$; $Z$ determines a unique line
$L\not\subset SecX$.  If $L$ intersects $X$ in a single point, then 
the restriction of $\varphi$ to $L$ resolves to a linear embedding
of $L$.  If $L$ does not intersect $X$, then the surjectivity of the
restriction map implies that $\varphi$ is an embedding along $L$.
Hence points and tangents are separated off $SecX$.
\end{proof}

In other words, we need only avoid the case where, after a choice 
of coordinates on $L$, the restriction of $\varphi$ to $L$ is the
system generated by $x^2,y^2$.  We introduce a 
condition that guarantees this does not happen.

\begin{defn}\label{kd}
A subscheme $X\subset\P^n$ {\bf \boldmath satisfies condition 
$(K_d)$} if $X$ is scheme theoretically cut out by forms 
$F_0,\ldots ,F_s$ of degree $d$ such that the trivial (or Koszul)
relations among the $F_i$ are generated by linear syzygies.
\end{defn}

More generally, let $V\subseteq H^0(\O_{\P^n}(d))$ be a linear system
of forms of degree $d$ with (possibly empty) base scheme $X$.
Then the pair $(X,V)$ {\bf \boldmath satisfies condition $(K_d)$} if the
trivial relations among the elements of $V$ are generated by
linear syzygies.  Write $(X,F_i)$ for the
pair $(X,V)$ if the set $\{ F_i\} $ generates the linear system $V$.
Perhaps the simplest example of varieties that do {\em not} satisfy $(K_d)$
is that of scheme theoretic complete intersections of hypersurfaces of
degree $d$.

\begin{rem}{Remark}{n2vsk2}
For a projective variety $X\subset\P^n$, M. Green \cite{mgreen} defines 
{\bf \boldmath condition $(N_2)$} as: $X$ projectively normal, 
ideal theoretically defined by quadrics $F_i$, and 
all of the syzygies among the $F_i$ are generated by linear ones.  
Examples include a smooth curve embedded by a line bundle
of degree at least $2g+3$ \cite{mgreen}, canonical curves with
Clifford index at least $3$ \cite{fsch},\cite{v}, Veronese embeddings
of $\P^n$, e.g. \cite{einlaz}, and all sufficiently large embeddings
of any projective variety \cite{mgreen},\cite{inamdar}.  

Clearly, if $X$ satisfies $(N_2)$, then $X$ satisfies the weaker
condition $(K_2)$.  Though $(K_2)$ is 
a technically simpler condition and will arise naturally, in practice most
examples we consider that satisfy $(K_2)$ will actually 
satisfy the stronger, and well studied, condition $(N_2)$.  As such,
we have made no attempt to understand examples where $X$ satisfies
$(K_2)$ but not $(N_2)$.  It is not difficult, however, to see where
such examples might arise.  Specifically, if $X\subset\P^n$ is a
smooth surface satisfying $(N_2)$ and if $h^1(\O_X(3))<h^1(\O_X(1))$,
then a general quadric section of $X$ is not projectively normal, but
is certainly scheme theoretically defined by quadrics.  

In light of Lemma~\ref{linearsubspace} below, one may also expect to find
examples by taking hyperplane sections of non arithmetically
Cohen-Macaulay varieties.
\end{rem}  
An advantage of the weaker 
condition $(K_d)$ is the following simple:
\begin{lemma}\label{linearsubspace}
Let $V$ be a linear system on $\P^n$ that satisfies $(K_d)$, and let 
$M\cong \P^k$ be a linear subspace of $\P^n$.  Then the restriction of 
$V$ to $M$ satisfies $(K_d)$.
{\nopagebreak \hfill $\Box$ \par \medskip}
\end{lemma}

This gives
\begin{cor}\label{weakembedding}
Let $X\subset\P^n$ be scheme theoretically defined by 
quadrics $F_0,\ldots ,F_s$ satisfying $(K_2)$.  Then the induced 
map $\varphi$ is an embedding off $SecX$.
\end{cor}

\begin{proof}
Let $L\subset\P^n$, $L\not\subset SecX$, be a line not intersecting 
$X$.  By Lemma~\ref{linearsubspace}, the restriction of the $F_i$
to $L$ must satisfy $(K_2)$.  However, it is easy to check that the 
only base point free system of quadrics on $\P^1$ satisfying $(K_2)$
is the complete system of quadrics.  Hence by 
Proposition~\ref{embeddingequiv}, $\varphi$ is an embedding off $SecX$.
\end{proof}

\begin{rem}{Remark}{hulek-oxbury}
A similar result was discovered independently by K. Hulek and
W. Oxbury \cite{ho}.
{\nopagebreak \hfill $\Box$ \par \medskip}
\end{rem}

Recall that a base point free linear system $W$ on $X$ is said to be {\em
$k$-very ample} if every zero dimensional subscheme of $X$ of length
$k$ spans a $\P^{k-1}$ in $\P(W)$.
We record the following elementary:
\begin{lemma}\label{4sep}
Assume $X\subset \P(W)=\P^n$ satisfies condition $(K_2)$ 
and contains no lines or conics.  
Then $W$ is a $4$-very ample linear system on $X$.
\end{lemma}

\begin{proof}
Assume to the contrary that there is a $2$-plane $M$ that intersects $X$ in
a scheme $Z$ of length $k\geq 4$.  Note that by
hypothesis $M$ cannot intersect $X$ in a scheme of positive dimension.

Two conics in $M$ intersect in a scheme of length $4$ if and only if they
have no common component.  However, a pair of plane conics cannot
satisfy $(K_2)$ unless they share a linear factor, which would imply a
positive dimensional base scheme.  Therefore, there can be no such $2$-plane.
\end{proof}

We return to the case of arbitrary $d$.
Via the closure of the graph
$\overline{\Gamma_{\varphi}}\subset \P^n\times \P^s$, we have a resolution
of the rational map $\varphi$:
\begin{center}
{\begin{minipage}{1.5in}
\diagram
  \overline{\Gamma_{\varphi}} \dto_{\pi_1} \drto^{\wtv} & \\
  \P^n  \rdashed|>\tip^{\varphi}  & \P^s
\enddiagram
\end{minipage}}
\end{center}
where $\wtv$ is the restriction of the projection onto the second factor. 
Note that $\wtv$ can be identified with the morphism
$\blow{\P^n}{X}\rightarrow \P^s$ induced by lifting the appropriate
sections of $\Gamma(\P^n,\I_X(d))$ to sections of 
$\Gamma(\wtp,\O_{\wtp}(dH-E))$.
\begin{prop}\label{fibers}
Let $(X,F_i)$ satisfy $(K_d)$ and assume that $X$ does not contain a
line.  If $a\in Im~{\wtv}$, then ${\wtv}^{-1}(a)=\P^k\times \{ a\}
  \subset \P^n \times \P^s$, where either
\begin{enumerate}
\item $k=0$ or 
\item $\pi_1(\P^k\times \{ a\} )=\P^k\subset \P^n$ intersects $X$ in a
hypersurface of degree $d$ in $\P^k$.
\end{enumerate}
\end{prop}

\begin{proof}
Take coordinates $[z_0,\ldots,z_n;t_0,\ldots,t_s]$ on $\P^n\times \P^s$
and let $S$ be the scheme defined by the equations
$\left\{ F_it_j-F_jt_i=0 \right\} \forall i,j$.  Clearly, 
$\overline{\Gamma_{\varphi}}\subset S$ as schemes and it is easy to verify
that $\overline{\Gamma_{\varphi}}=S$ off of $E$, the exceptional divisor of
the blow up.  

Now, considering a syzygy as a vector of forms, 
let $$\displaystyle \left\{ (a_{\ell 0},\ldots
,a_{\ell s}) ; 0\leq \ell \leq r\right\} $$ 
generate the linear syzygies among the $F_i$. Let 
$T\subset \P^n\times \P^s$ be the subscheme defined by the 
equations $\displaystyle \left\{ \sum_{k=0}^s a_{\ell
k}(z)t_k ; 0\leq \ell \leq r\right\} $.  Again it is clear that
$\overline{\Gamma_{\varphi}}\subset T$ as schemes, 
(Cf. \cite[\S 1]{hks}) but condition
$(K_d)$ implies that $T\subset S$, hence 
$\overline{\Gamma_{\varphi}}=T$ off of $E$.

$\wtv^{-1}(a)$ is contained as a scheme in
$T_a$, where $T_a$ is the fiber over $a$ of the projection map
restricted to $T$.  Without loss of generality, assume
$a=[1,0,\ldots ,0]\in \P^s$.  $T_a$
is then scheme theoretically defined by the bihomogeneous equations:
$$\displaystyle \left\{ \sum_{k=0}^s a_{0k}t_k ,\ldots ,\sum_{k=0}^s
a_{rk}t_k ,t_1,t_2,\ldots ,t_s \right\} $$
and so more simply by:
$$\displaystyle \left\{ a_{00}t_0,\ldots ,a_{r0}t_0,t_1,t_2,\ldots
,t_s\right\} $$
giving:
\begin{eqnarray*}
T_a & = & \P^k \times [1,0,\ldots ,0]\subset \P^n \times \P^s \\
 & \cong & \P^k 
\end{eqnarray*}
where $\P^k$ is the linear subspace of $\P^n$ defined by the $\{ a_{\ell
0}\} $.  We have just seen $\wtv^{-1}(a)=T_a$ off of $E$;
$T_a$ is irreducible, however, so $\wtv^{-1}(a) = T_a \cong \P^k$
as long as either of the following is true: 
\begin{enumerate}
\item $T_a$ is a reduced point (i.e. $k=0$).
\item $T_a$ and $\wtv^{-1}(a)$ are not both contained in
$E$.
\end{enumerate}

To guarantee that the second possibility occurs if $T_a$ has
positive dimension, note that $X$ does not contain a line.  Therefore
$\pi_1(T_a)$, which is a linear subspace of $\P^n$, cannot be isomorphic
to a positive-dimensional reduced linear subscheme of $X$.  Hence $T_a$
cannot be contained in $E$. 
\end{proof}

\begin{rem}{Remark}{alternatehyp}
Rather than the hypothesis that $X$ contain no lines, one could instead
insist that $X$ is not set theoretically cut out by any subsystem
of the $F_i$.  We choose the former hypothesis as it will be
necessary below.
{\nopagebreak \hfill $\Box$ \par \medskip}
\end{rem}

Denoting by $Sec^1_dX$ the variety of lines intersecting $X$ in a 
subscheme of length at least $d$,
we have a natural extension of Corollary~\ref{weakembedding}:

\begin{thm}\label{embed}
Let $(X,F_i)$ be a pair that satisfies $(K_d)$, and assume that $X$ does
not contain a line.  Then $\wtv$ is an embedding off
the proper transform of $Sec^1_dX$.
Furthermore, if $X$ is smooth then the image of $\wtv$
is a normal subvariety of $\P^s$. 
\end{thm}

\begin{proof}
To prove the first claim, note that because the fibers of
$\wtv$ are reduced, we need only show that points are
separated. 

Take $p,q\in \overline{\Gamma_{\varphi }}$, and assume
$\wtv(p)= \wtv(q)=\{ r\} \in\P^s$.  Then
by Proposition~\ref{fibers}, $S=\wtv^{-1}(r)$ satisfies
$\pi(S)=\P^k, k>0$.  There are then two possibilities:
\begin{enumerate}
\item $\pi(S)\cap X$ is a $d$-ic hypersurface in $\pi(S)$, and hence
every line in $\pi(S)$ is a $d$-secant line of $X$, which implies that
$\pi(S)\subseteq Sec^1_d(X)$, and so $p$ and $q$ are in the proper
transform of the variety of $d$-secant lines.

\item $\pi(S)\subseteq X$.  In this case,
however, $\pi(S)$ is a positive dimensional linear subvariety of $X$,
which is not allowed by hypothesis.
\end{enumerate}

To see that the image is normal if $X$ is smooth, notice that
by identification with the blow up of $\P^n$ along $X$, 
$\overline{\Gamma_{\varphi }}$ is smooth, hence normal.
We have just shown that the fibers are reduced and connected, hence the
result follows from the fact that the image of a proper 
morphism from a normal variety with reduced, connected fibers is normal.
\end{proof}

\begin{rem}{Remark}{phianembedding}
The proof of Theorem~\ref{embed} implies that if $X\subset \P^n$ is
scheme theoretically defined by forms of degree $d$ that satisfy
$(K_d)$, then the map $\varphi :\P^n\setminus X\rightarrow \P^s$ is an
embedding off of $Sec^1_dX$, even if $X$ {\em does} contain a line,
extending Corollary~\ref{weakembedding}.
{\nopagebreak \hfill $\Box$ \par \medskip}
\end{rem}

Theorem~\ref{embed} provides a simple extension of a vanishing
result for powers of ideal sheaves of projective
varieties in \cite[1.10]{bel}.  The bound is only improved by
one degree; this comes precisely from the fact that we know
$\O(dH-E)$ is big.
\begin{cor}\label{littlevanishing}
If $X$ is smooth, irreducible and satisfies $(K_d)$ then 
Theorem~\ref{embed} 
shows that the linear system $\O(dH-E)$ on 
$\overline{\Gamma_{\varphi }}=\blow{\P^n}{X}$
is big and nef, and $\O(kH-mE)$ is very ample if $\frac{k}{m}>d$.  
A simple application 
of the Kawamata-Viehweg vanishing theorem gives
$$H^i\left( \P^n,\I_X^a(k)\right) =0, i>0 , k\geq d(e+a-1)-(n+1)$$
where $e$ is the codimension of $X$ in $\P^n$.
{\nopagebreak \hfill $\Box$ \par \medskip}
\end{cor}

\section{Results on Secant Varieties}\label{secresults}

We describe a vector bundle on $\H^{2}X=Hilb^{2}(X)$ and a morphism
to projective space giving rise to $SecX$.  Our construction 
follows that of \cite[\S 1]{bertram1}, \cite[VIII.2]{acgh}, and \cite{sch}
where this is done for curves with the identification of $\H^{2}X$
with $S^{2}X$.

Let $V\subseteq \Gamma(X,L)$ be a very ample linear system and
denote by $\D$ the universal
subscheme of $X\times \H^{2}X$ and note $\D\cong \blow{X\times X}{\Delta}$.
Let $\pi:X\times \H^{2}X\rightarrow
X$ and $\pi_{2}:X\times \H^{2}X\rightarrow \H^{2}X$ be the
projections, and let $L$ be any line bundle on $X$.  Form the
invertible sheaf $\O_{\D}\tensor \pi^*L$ on $\D
\subset X\times \H^{2}X$. Now $\pi_{2}|_{\D}:\D
\rightarrow \H^{2}X$ is flat of degree $2$, hence
$\mathscr{E}_L=(\pi_{2})_*(\O_{\D}\tensor \pi^*L)$ is a
locally free sheaf of rank $2$ on $\H^{2}X$.  We define the {\bf first 
secant bundle} of $X$ with respect to $L$ to be the $\P^1$-bundle
$B^1(L)=\P_{\H^{2}X}(\mathscr{E}_L)$. 

To define the desired map, push the natural restriction  
$\pi^*L\rightarrow \O_{\D}\tensor \pi^*L$
down to $\H^{2}X$ giving an evaluation map
$H^0(X,L)\tensor \O_{\H^{2}X}\rightarrow \mathscr{E}_L$
which in turn for any linear system $V\subseteq H^0(X,L)$ restricts to
$V\tensor \O_{\H^{2}X}\rightarrow \mathscr{E}_L$.
Now a fiber of $\mathscr{E}_L$ over a point $Z\in \H^{2}X$ is
$H^0(X,L\tensor \O_Z)$, so if $V$ is very ample then this map
is surjective and we obtain a morphism:
$$\beta_1:B^1(L)\rightarrow \P(V)\times \H^{2}X\rightarrow \P(V)$$
The image of this morphism
is the secant variety to $X$ in $\P(V)$. 

\begin{rem}{Remark}{intograss}
It will be useful to note that the above surjection also induces a morphism 
$\H^2X\rightarrow\G(1,V)$ which is an embedding as long as $V$ is 
$3$-very ample \cite{cat-gott}.   
{\nopagebreak \hfill $\Box$ \par \medskip}
\end{rem}

\begin{rem}{Notation}{switchtosigma}
To help simplify notational clutter, we denote $\Sigma=Sec^1_2X$ and
$\mathcal{T}=TanX$.  This should cause no confusion as we will be concerned
with a fixed variety $X$, and will be primarily concerned only with
the first secant variety.
{\nopagebreak \hfill $\Box$ \par \medskip}
\end{rem}

\begin{rem}{Hypothesis}{hyp}
For the remainder of this section, $X\subset \P^n$  will denote a
smooth, irreducible, non-degenerate variety, scheme theoretically
defined by quadrics $(F_0,\ldots ,F_s)=V\subset\Gamma(\P^n,\O(2))$ 
satisfying $(K_2)$.  Assume that $X$
{\it contains no lines and no conics}.  In particular,
the embedding of $X$ is $4$-very ample.
{\nopagebreak \hfill $\Box$ \par \medskip}
\end{rem}

In this situation, Theorem~\ref{embed} implies that the map $\wtv$ is an
embedding off the proper transform of the secant variety to $X$.
Here we study what the map does when restricted to the proper transform
of the secant variety.  Denote by
$\wts$ the proper transform of the secant 
variety under the blowing up $\pi:\wtp\rightarrow \P^n$ of
$\P^n$ along $X$.  By a slight abuse of notation,
write $\wtv: \wts \rightarrow \P^s$ for the restriction of
$\wtv: \wtp\rightarrow \P^s$. 
We show first that there is a map $\wts \rightarrow \H^2X$, and then an
embedding of $\H^2X$ into the image variety of $\wtv$, such that the
composition factors $\wtv :\wts\rightarrow \P^s$ and the map
$\wts \rightarrow \H^2X$ is a $\P^1$-bundle.

\begin{rem}{Remark}{P1s}
Because of the assumption that $X$ contains no lines and no conics, 
each fiber of $\wtv: \wts
\rightarrow \P^s$ is isomorphic to $\P^1$ by Proposition~\ref{fibers}.
In particular, given a point in $\wts$ or in $\SecX\setminus X$, one
can say on which secant or tangent line it lies. 
{\nopagebreak \hfill $\Box$ \par \medskip}
\end{rem}

\begin{lemma}
There is a morphism $g:\wts\rightarrow \H^2X$ taking a point $p$
to the length $2$ subscheme $Z$ of $X$ 
determining the secant line on which $p$ lies.
\end{lemma}

\begin{proof}
We construct a morphism 
$\wts\rightarrow\G(1,n)$ whose image is $\H^2X$.
Let $Y=Im~\wtv$, and push the surjection
$$H^0(\wts,\O_{\wts}(H))\tensor \O_{\wts}\rightarrow 
\O_{\wts}(H)\rightarrow 0$$
down to $Y$:
$$H^0(\wts,\O_{\wts}(H))\tensor \O_Y\rightarrow \wtv_*\O_{\wts}(H)$$
The sheaf $\wtv_*\O_{\wts}(H)$ is locally free of rank $2$, and
the map is surjective as $\O(H)$ maps a fiber of $\wtv$ to a linearly
embedded $\P^1\subset \P^n$.  Pulling this surjection back to $\wts$
gives a surjection from a free rank $n+1$ sheaf to a rank $2$ vector bundle, 
hence a morphism $\wts\rightarrow \G(1,n)$ taking a fiber of $\wtv$
to the point representing the associated secant line.  The image of
this morphism is clearly $\H^2X\hookrightarrow \G(1,n)$
from Remark~\ref{intograss}.
\end{proof}

Having constructed a map $g:\wts\rightarrow \H^2X$, we 
construct an embedding $f:\H^2X\hookrightarrow \P^s$ so that
$\wtv=(f\of g):\wts\rightarrow \P^s$. 

Let $Z\in \H^2X$ be a length $2$ subscheme of $X$, and let
$\ell_Z\subset \P^n$ be the line determined by $Z$.  Note that by
hypothesis $\ell_Z$ does not lie on $X$.  There are homomorphisms:
$$r_Z:V\rightarrow H^0(\P^n,\O_{\P^n}(2)\tensor \O_{\ell_Z})$$
$f$ is set theoretically given by associating to every $Z\in \H^2X$ the 
$1$-dimensional quotient $\gmodh{$V$}{$ker(r_Z)$}$.  
This association is injective by Remark~\ref{P1s}.

\begin{lemma}\label{hilbintoPV}
$f:\H^2X\rightarrow \P(V)=\P^s$ is a morphism.
\end{lemma}

\begin{proof}
Let $L=\O_{\wtp}(2H)$ and form $\pi_1^*L$ on
$\wtp\times \H^2X$.  Embed:
\begin{eqnarray*}
\wts & \hookrightarrow & \wtp\times \H^2X\\
p & \mapsto & (i(p),g(p))
\end{eqnarray*}
where $i:\wts\hookrightarrow \wtp$ is the inclusion.  
Applying $\pi_{2_*}$ to the surjection
$\pi_1^*L\rightarrow \pi_1^*L\tensor \O_{\wts}\rightarrow 0$  
gives a map:
$$H^0(\wtp,L)\tensor \O_{\H^2X}\rightarrow
\pi_{2_*}\left( \pi_1^*L\tensor \O_{\wts}\right) $$
Recalling $V\subseteq H^0(\wtp,L\tensor \O_{\wtp}(-E))$, there
is a map: 
$$V\tensor \O_{\H^2X}\rightarrow \pi_{2_*}\left( \pi_1^*L\tensor
\O_{\wts}\right) $$
where a fiber of the coherent sheaf $\pi_{2_*}\left( \pi_1^*L\tensor
\O_{\wts}\right) $ over a point $Z\in \H^2X$ is isomorphic to
$H^0(\P^n,\O_{\P^n}(2)\tensor \O_{\ell_Z})$.  By the above remarks, this 
map has rank $1$, hence gives a surjection to a line bundle on $\H^2X$
with fiber over $Z$ isomorphic to $H^0(\P^n,\O_{\P^n}(2)\tensor \I_X
\tensor \O_{\ell_Z})$.  $f$ is the morphism induced by this surjection.
\end{proof}

The diagram
\begin{center}
{\begin{minipage}{1.5in}
\diagram
  \wts \rto^{g}
\drto^{\wtv} & \H^2X \dto^{f} \\
  & \P^s
\enddiagram
\end{minipage}}
\end{center}
then commutes, and because the fibers of
$\wtv$ are reduced, those of $f$ are as well:

\begin{prop}\label{hilbintoPs}
Let $(X,V)$ be a pair satisfying $(K_2)$, and assume $X\subset \P^n$
is smooth, irreducible, and contains no lines or conics.
Then the morphism $f:\H^2X\hookrightarrow \P(V)$ above is an
embedding. 
{\nopagebreak \hfill $\Box$ \par \medskip}
\end{prop}

\begin{lemma}\label{FisD}
With hypotheses as in Proposition~\ref{hilbintoPs}, the exceptional divisor
of the blow up $\wts\rightarrow \SecX$ is isomorphic to 
$\blow{X\times X}{\Delta}$.
\end{lemma}

\begin{proof}
Let $F\subset \wts$ be the exceptional divisor of this
blow up.  Let $Y$ be the image of:
\begin{eqnarray*}
 F & \rightarrow & X\times \H^2X \\
 p & \mapsto & \left( \pi (p),\wtv(p)\right)
\end{eqnarray*}
$Y$ is flat of degree $2$ over $\H^2X$.  Indeed, by the structure
of $\wtv$ the fiber over a point in $\H^2X$ is exactly the
corresponding length $2$ subscheme of $X$, hence
$Y$ induces the identity morphism
$id:\H^2X\rightarrow \H^2X$.  By the universal property of $\D\cong
\blow{X\times X}{\Delta}$, the universal subscheme of $X\times \H^2X$,
we have:
$$Y\cong (id_X\times id_{\H^2X})^{-1}(\D)$$
The map from $F$ to $Y$ is
a finite birational morphism to a smooth variety, so is an
isomorphism, hence $F\cong \blow{X\times X}{\Delta}$.
\end{proof}

This allows
another construction of the secant bundle $B^1(L)$: Writing
$\picard\wtp\cong \Z H+\Z E$, form the line bundle $H\tensor
\O_{F}$ on the exceptional divisor $F\subset \wts$.  The
restriction of $\wtv$ to $F$ is a degree two map to
$\H^2X$.  Let $\E=\wtv_*(H\tensor \O_{F})$.  By the
identification of $F$ with $\D$ and of $H\cong \pi^*\O_{\P^n}(1)$ with
$L$ on $X$, we see $B^1(L)\cong \P_{\H^2X}(\E)$.

\begin{thm}\label{getabundle}
Let $X\subset \P^n$ be smooth, irreducible and satisfy $(K_2)$.  If
$X$ contains no lines and no conics then
$\wtv:\wts\rightarrow\H^2X$ is the $\P^1$-bundle 
$\P_{\H^2X}(\E)\rightarrow \H^2X$.
\end{thm}

\begin{proof}
Note that $\E=\wtv_*(H\tensor \O_{F})$ and
$\wtv_*(H\tensor \O_{\wts})$ are isomorphic
rank two vector bundles on $\H^2X$ and that $H\tensor
\O_{\wts}$ is generated by its global sections.  Hence 
there is a surjection $\wtv^*\E\rightarrow H\tensor
\O_{\wts}\rightarrow 0$ which induces a morphism 
$\kappa :\wts\rightarrow \P_{\H^2X}(\E)$.

This gives the diagrams:
\begin{center}
{\begin{minipage}{1.5in}
\diagram
  \wts \rto^{\kappa} \dto^{\pi} & \P_{\H^2X}(\E) \dlto^{\beta_1} &
\wts \rto^{\kappa} \drto^{\wtv} & \P_{\H^2X}(\E)
\dto^p\\
  \SecX & & & \H^2X
\enddiagram
\end{minipage}}
\end{center}
where $p$ is the natural projection map.  $\kappa$ makes both triangles
commute, and so is a finite (by the second diagram) birational (by the
first) morphism to a smooth variety, hence an isomorphism.
\end{proof}

\begin{cor}\label{alliswell}
Under the hypotheses of Theorem~\ref{getabundle}:
\begin{enumerate}
\item $\wts$ is smooth and $\SecX$ is smooth off $X$
\item $\wt{\TanX}$ is smooth and $\TanX$ is smooth off $X$
\item $\SecX$ is normal
\end{enumerate}
\end{cor}

\begin{proof}
The smoothness of $\wts$ is immediate from Theorem~\ref{getabundle}.

For the second claim, note
that $\wtv$ maps $\wt{\TanX}$ to the diagonal in
$\H^2X$, which is the projectivized tangent bundle to $X$, hence
smooth (by the {\em diagonal} in $\H^2X$, we mean the proper transform
on the diagonal under the birational morphism $\H^2X\rightarrow S^2X$).
Therefore, $\wt{\TanX}$ is smooth as above. 

To show $\SecX$ is normal, we use the (just proven) fact that
$\wts$ is smooth, hence normal.  It suffices
to show that for $p\in X$, $\pi^{-1}(p)$ is reduced and connected where
$\pi :\wts\rightarrow \SecX$ is the blow up of $\SecX$
along $X$.  But by Lemma~\ref{FisD}, $\pi^{-1}(p)\cong \blow{X}{p}$.
\end{proof}

\begin{rem}{Remark}{workofbertram}
In \cite{bertram1}, A. Bertram shows directly that
$\wts$ is isomorphic to $B^1(L)$ if $X$ is a smooth curve
embedded by the complete linear system associated to a line bundle $L$
that is $4$-very ample.
{\nopagebreak \hfill $\Box$ \par \medskip}
\end{rem}

\begin{rem}{Remark}{deficientsecants}
A simple consequence of the above results is that a smooth, irreducible
variety $X\subset \P^n$ satisfying $(K_2)$ with no lines and no conics 
has a non-deficient secant variety.  It can 
be shown \cite[3.6.1]{vermeire} that if $\delta =2r+1-dim(\SecX)$ 
is the deficiency of the
secant variety to $X$, then
$\delta =\frac{2r-dim Y}{2}$
where $Y$ is the image variety of $\wtv:\wts\rightarrow\P^s$.

In particular, this shows that the dimension of $\SecX$ is determined
by the dimension of the fibers of $\wtv$: The generic pair of points
of $X$ lies on a quadric hypersurface of (maximal) dimension 
$d$ if and only if $dim(\SecX)=2r+1-d$.
{\nopagebreak \hfill $\Box$ \par \medskip}
\end{rem}

\section{Geometric Flip Construction}\label{flipconst}

In this section we are motivated by \cite{thad} to construct
a flip centered about $\H^2X$.  We recall his construction.

\subsection{Work of Thaddeus}

In \cite{thad} Thaddeus considers the moduli problem of semi-stable 
pairs $(E,s)$ consisting of a rank 
two bundle $E$ with $\wedge^2E=\Lambda$, and a section 
$s\in \Gamma(X,E)-\{ 0\}$.  This is interpreted as a GIT problem, and
by varying the linearization of the group action, a collection
of (smooth) moduli spaces $M_1,M_2,\ldots,M_k$, $k=\left
[ \frac{d-1}{2}\right] $, 
is constructed.  As stability is an open condition, these spaces
are birational.  In fact, they are isomorphic in codimension one, and
may be linked via a diagram
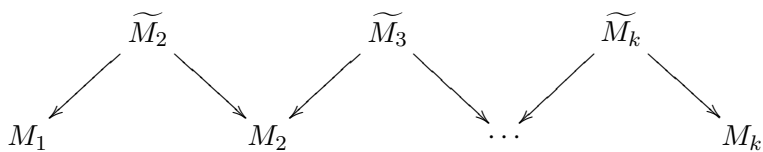
\begin{center}
{\begin{minipage}{1.5in}
\diagram
 & \wt{M_2} \dlto \drto & & \wt{M_3} \dlto \drto & & \wt{M_k} \dlto \drto & \\
 M_1  & & M_2 & & \cdots & & M_k
\enddiagram
\end{minipage}}
\end{center}
where there is a morphism $M_k\rightarrow M(2,\Lambda)$.
The relevant observations are first that this is a diagram of flips (in fact
it is shown in \cite{bertram4} that it is a sequence of log flips)
where the ample cone of each $M_i$ is known,
second that $M_1$ is the blow up of $\P(\Gamma(X,K_X\tensor\Lambda)^*)$
along $X$, and finally that $\wt{M_2}$ is the blow up of $M_1$ along
the proper transform of the secant variety and that all of the
flips can be seen as blowing up and down various higher secant varieties.

Our inspiration can be stated as follows:  The sequence of flips
in Thaddeus' construction, constructed via Geometric Invariant Theory, 
can be realized as a sequence of natural geometric constructions depending
only on the original embedding of $X\subset \P^n$.  An
advantage of this approach is that the smooth curve $X$ can
be replaced by any smooth variety.  Even in the curve case, ours 
applies to situations where the original construction does not hold 
(e.g. for canonical curves with Clifford index at least $3$).

Thaddeus goes on \cite[7.8]{thad} to compute the dimension of
the spaces $H^0(\P H^0(K_X\tensor \Lambda),\I^a_X(k))$ for
certain values of $d,g,a,k$.  In particular, this computation is
used to verify the rank-two Verlinde formula.  A part
of our motivation is to try to extend this computation
to a much larger class of varieties.

\subsection{Outline of Our Construction}
With notation as above assume
$X\subset \P^n$ is smooth, irreducible, satisfies $(K_2)$, 
and contains no lines or conics.  
Let $r=dim~X$ and assume that $n-2r-1\geq 2$,
i.e. that $\SecX$ is not a hypersurface in $\P^n$.  Write $\P(\E)$
for the secant bundle $\P_{\H^2X}(\E)$ and identify $\wts$ with $\P(\E)$
(Theorem~\ref{getabundle}).  Simply write
$\wtv:\P(\E)\rightarrow \H^2X$ for the restriction of
$\wtv:\wtp\rightarrow \P^s$.

We begin with a general construction:
Let $f:X\dashrightarrow Y$ and $g:X\dashrightarrow Z$ be rational maps 
of irreducible varieties,
$f$ birational.  Form $\overline{\Gamma_{f,g}}\subset X\times Y\times Z$, 
the closure of the graph of $(f,g):X\dashrightarrow Y\times Z$, and
let $M_2\subset Y\times Z$ be the image of $\overline{\Gamma_{f,g}}$
under the obvious projection.  We have
the following diagrams where all maps are projections, the left
included in the right by restriction:
\begin{center}
{\begin{minipage}{1.5in}
\diagram
 & \overline{\Gamma_{f,g}} \dlto \drto & & & X\times Y\times Z \dlto \drto & \\
 \overline{\Gamma_f} \drto & & M_2 \dlto & X\times Y \drto & & 
Y\times Z \dlto \\
 & Y & & & Y &
\enddiagram
\end{minipage}}
\end{center}
Note that $X, Y, \overline{\Gamma_f}, \overline{\Gamma_{f,g}}$, and $M_2$ 
are all birational.

Of particular interest is the case where 
$\overline{\Gamma_f}\rightarrow Y$ is a small morphism with
exceptional locus $W$.  Assume there exists a line bundle 
$L$ on $\overline{\Gamma_f}$ with base scheme $W$ and take
$Z=\P H^0\left( \overline{\Gamma_f},L\right) $.  Then 
$\overline{\Gamma_f}$ and $M_2$ are isomorphic in codimension
one.  Furthermore, if $\overline{\Gamma_f}$ and $M_2$ are
factorial varieties, then the
image of $L$ under the isomorphism $\picard \overline{\Gamma_f}\cong
\picard M_2$ is a globally generated line bundle on $M_2$.

We give an explicit construction of this
birational transformation: take 
$\overline{\Gamma_f}\rightarrow Y$ above 
to be $\wtv:\wtp\rightarrow \P^s$.  As the exceptional locus of $\wtv$ is
$\wts$, we find a line bundle on $\wtp$ 
whose base scheme is $\wts$, identify
explicitly the exceptional loci in the diagram, give
explicit linear systems defining the morphisms and finally show
that the space $M_2$ is smooth.

\subsection{The Diagram of Exceptional Loci}

An examination of Thaddeus' construction \cite[3.11]{thad} suggests
we identify a vector bundle $\F$ on $\H^2X$ of rank
$n-2r-1=codim(\P(\E),\wtp)$ such that:
\begin{eqnarray}
\wtv^*(\F) & \cong & N^*_{\P(\E)/\wtp} \tensor \O_{\P(\E)}(-1) \label{wantthis}
\end{eqnarray}
and then construct one of the exceptional loci as $\P_{\H^2X}(\F)$.
Writing $N_{\P(\E)}(k)=N_{\P(\E)/\wtp}\tensor \O_{\P(\E)}(k)$, we verify
$\F =\wtv_*N^*_{\P(\E)}(-1)$ satisfies (\ref{wantthis}).

\begin{prop}\label{ample}
Let $\wt{Y}\cong \P^1$ be a fiber of
$\wtv:\P(\E)\rightarrow \H^2X$.  Then $N_{\P(\E)}(k) \tensor \O_{\wt{Y}}\cong 
\oplus \O_{\wt{Y}}(k-1)$ and it follows that
$$\wtv^*\wtv_*N^*_{\P(\E)}(-1)\cong N^*_{\P(\E)}(-1)$$
\end{prop}

\begin{proof}

Denoting tangent sheaves by $\T$, it is easy to see that
$$(\pi^*\T_{\P^n})\tensor\O_{\wt{Y}}\cong \T_{\P^n}\tensor\O_Y
\cong \O(2)\oplus\O(1)\oplus\cdots\oplus\O(1)$$ Let $F$ be the 
universal quotient bundle on the exceptional divisor $E$ and 
$j:E\hookrightarrow \wtp$ the inclusion.
Because the intersection of $\wt{Y}$ with $E$ is a scheme of length two, 
pulling the exact sequence \cite[15.4]{fulton} 
$\ses{\T_{\wtp}}{\pi^*\T_{\P^n}}{j_*F}$
back to $\wt{Y}$ gives $\T_{\wtp}\tensor\O_{\wt{Y}}(k)$ ample
for $k\geq 2$.  The sequence
$$\ses{\T_{\P(\E)}\tensor\O_{\wt{Y}}}
{\T_{\wtp}\tensor\O_{\wt{Y}}}
{N_{\P(\E)/\wtp}\tensor\O_{\wt{Y}}}$$
then gives $N_{\P(\E)}(k)\tensor \O_{\wt{Y}}$ ample for $k\geq 2$
because it is the quotient of an ample bundle \cite[III.1.7]{hart2};
hence $N_{\P(\E)}\tensor \O_{\wt{Y}}\cong
\oplus \O_{\wt{Y}}(a_i)$ where the $a_i\geq -1$.  A 
straightforward computation of the determinant via the isomorphism 
$\omega_{\P(\E)}\cong \omega_{\wtp}\tensor\Lambda^{n-2r-1}N_{\P(\E)}$
shows $a_i=-1$, and the first part of the statement holds. 
The second follows immediately.
\end{proof}

We construct the diagram of exceptional loci for the flip:
Let $$E_2'=\P_{\P(\E)}\left( N^*_{\P(\E)}(-1)\right)
=\P_{\P(\E)}({\wtv}^*\F)$$ 
and $$E_2=\P_{\P(\E)}\left( N^*_{\P(\E)}\right) $$
Hence $E_2$ is the exceptional divisor of the blow up
$\blow{\wtp}{\P(\E)}$.  As the vector bundles
defining $E_2$ and $E_2'$ differ by the twist of a line
bundle, there is an isomorphism $\gamma $ \cite[II.7.9]{hart}: 
\begin{center}
{\begin{minipage}{1.5in}
\diagram
 E_2 \rto^{\gamma} \drto_{\pi} & E_2' \dto^{\pi'} \\
 & \P(\E)
\enddiagram
\end{minipage}}
\end{center}
with the property that
\begin{equation}\label{bookkeeping} 
\displaystyle \gamma^*(\O_{E'_2}(1)) \cong \O_{E_2}(1)\tensor 
\pi^*(\O_{\P(\E)}(-1))
\end{equation}

Writing $\P(\F)=\P_{\H^2X}(\F)$, there is a morphism 
$E_2'\rightarrow\P(\F)$ induced by the natural surjection $({\wtv}\of
\pi')^*(\F)\rightarrow \O_{E_2'}(1)\rightarrow 0$.  Via $\gamma $ we get a
morphism $h:E_2\rightarrow \P(\F)$ induced by the surjection 
(note (\ref{bookkeeping}))
$$(\wtv\of \pi)^*(\F) \rightarrow \O_{E_2}(1)\tensor
\pi^*\left( \O_{\P(\E)}(-1)\right) \rightarrow 0$$
and hence a diagram:
\begin{center}
{\begin{minipage}{1.5in}
\diagram
 & E_2 \dlto_{\pi} \drto^{h} & \\
  \P(\E) \drto_{\wtv} & & \P(\F) \dlto^f \\
 & \H^2X &
\enddiagram
\end{minipage}}
\end{center}
The isomorphism $\gamma$ gives $E_2\cong \P_{\P(\E)}(\wtv^*\F)$.  Note
the following symmetry property:
\begin{lemma}\label{alsoabundle}
$E_2\cong \P_{\P(\F)}(f^*\E)$.
\end{lemma}

\begin{proof}
To give a map $E_2\rightarrow \P_{\P(\F)}(f^*\E)$ it is
equivalent to give a surjection $h^*f^*\E\rightarrow
\mathscr{K}\rightarrow 0$ for some line bundle $\mathscr{K}$ on $E_2$.
By the above diagram, this is equivalent to a surjection
$\pi^*\wtv^*\E\rightarrow \mathscr{K}\rightarrow
0$, which we obtain from the natural surjection
$\wtv^*\E\rightarrow \O_{\P(\E)}(1)\rightarrow 0$ on
$\P(\E)$.  As the fibers of $h$ are isomorphic to $\P^1$, it is clear
that the induced map is an isomorphism.
\end{proof}

Let the very ample invertible sheaf $M$
on $\H^2X\subset \P^s$ be the restriction of $\O_{\P^s}(1)$.  Then for
every $k$ sufficiently large, $\O_{\P(\F)}(1)\tensor f^*M^k$ is very
ample on $\P(\F)$, \cite{hart},Ex. II.7.14], and so gives an embedding
$i:\P(\F)\hookrightarrow \P^r$.  The induced morphism $i\of
h:E_2\rightarrow \P^r$ is given by a linear system associated to
the line bundle: 
\begin{eqnarray}
(i\of h)^*\left( \O_{\P^r}(1)\right) & \cong & h^*\left(
\O_{\P(\F)}(1)\tensor f^*M^k\right) \nonumber\\
& \cong & h^*f^*M^k\tensor \O_{E_2}(1)\tensor \pi^*\left(
\O_{\P(\E)}(-1)\right) \label{pullback}  
\end{eqnarray}
Since $h_{*}\O_{E_2}=\O_{\P(\F)}$ by Lemma~\ref{alsoabundle}, the 
projection formula yields: 
$$\Gamma\left( \P(\F), \O_{\P(\F)}(1)\tensor f^*M^k\right) =
\Gamma\left( E_2,h^*\left( \O_{\P(\F)}(1)\tensor f^*M^k\right) \right) $$
hence:
\begin{lemma}\label{itscomplete}
The complete linear system $|\O_{\P(\F)}(1)\tensor f^*M^k|$ on $\P(\F)$ pulls
back to the {\bf \em complete} linear system
$|h^*(\O_{\P(\F)}(1)\tensor f^*M^k)|$ on $E_2$.
{\nopagebreak \hfill $\Box$ \par \medskip}
\end{lemma}

\subsection{The Total Spaces}
We build the total spaces of the flip containing the
diagram of exceptional loci, with those maps given by restriction.
Three of the four spaces have been constructed
already: $\wtp,Im~\wtv, $ and $\blow{\wtp}{\P(\E)}$.  We construct the
fourth (and most interesting!) as the image of a linear system on
$\blow{\wtp}{\P(\E)}$.  This construction proceeds in several steps:
First, we identify (\ref{rest}) an invertible sheaf on
$\blow{\wtp}{\P(\E)}$ that restricts to $h^*\left( \O_{\P(\F)}(1)\tensor
  f^*M^k\right)$ on $E_2$ (Cf. Lemma~\ref{itscomplete}).  We then show
that the associated complete linear system gives a 
birational morphism which is an embedding off $E_2$, and that
its restriction to $E_2$ is the complete linear system associated to
$h^*\left( \O_{\P(\F)}(1)\tensor f^*M^k\right)$.

Following the notation of \cite{thad}, denote 
$\wt{M_2}=\blow{\wtp}{\P(\E)}$.
Writing $$\picard\wt{M_2}=\Z H+ \Z E_1+ \Z E_2=\Z (\pi^*H)+\Z
(\pi^*E)+\Z E_2$$ 
and noting $\O_{E_2}(-E_2)=\O_{E_2}(1)$, we have:
\begin{eqnarray}
\O_{E_2}\left( (2k-1)H-kE_1-E_2\right) & \cong &
(f\of h)^*M^k\tensor \O_{E_2}(1)\tensor \pi^*\O_{\P(\E)}(-1)\nonumber \\
 & \cong & (i\of h)^*\left( \O_{\P^r}(1)\right) \label{rest}
\end{eqnarray}
by equation (\ref{pullback}).

\begin{prop}\label{free}
Let $\L$ be an invertible sheaf on a complete variety $X$, and 
let $\mathscr{B}$ be any locally free sheaf.
Assume that the map $\lambda:X\rightarrow Y$ 
induced by $\mathscr{L}$ is a birational morphism and that $\lambda$ 
is an isomorphism in a neighborhood of $p\in X$.
Then for all $n$ sufficiently large, the map
$$H^0(X,\mathscr{B}\tensor \mathscr{L}^n)\rightarrow
H^0(X,\mathscr{B}\tensor \mathscr{L}^n\tensor \O_p)$$
is surjective.
\end{prop}

\begin{proof}
Push the exact sequence
$$0\rightarrow \mathscr{B}\tensor \mathscr{L}^{n}\tensor \I_p\rightarrow
\mathscr{B}\tensor \mathscr{L}^{n}\rightarrow \mathscr{B}\tensor
\mathscr{L}^{n}\tensor \O_p\rightarrow 0$$
down to $Y$.  Because $\lambda$ is an isomorphism in a neighborhood of $p$,
the map $$\lambda_{*}(\mathscr{B}\tensor \mathscr{L}^{n}\tensor \O_p)\rightarrow 
R^1\lambda_{*}(\mathscr{B}\tensor \mathscr{L}^{n}\tensor \I_p)$$ 
is the zero map, hence there is an exact sequence on $Y$:
$$0\rightarrow  \lambda_{*}(\mathscr{B}\tensor
\I_p)(n)\rightarrow \lambda_{*}(\mathscr{B})(n)
\rightarrow  \lambda_{*}(\mathscr{B})(n)\tensor
\O_p\rightarrow 0$$
where $\mathscr{L}\simeq  \lambda^*\O_{Y}(1)$. 
Since $\O_{Y}(1)$ is (very) ample,
$H^1(Y,\lambda_{*}(\mathscr{B}\tensor \I_p)(n))=0$ for all $n$
sufficiently large.  Therefore there is a section of
$\lambda_{*}(\mathscr{B})(n)$ that does not vanish at $p\in Y$
which can be pulled back to a section of $\mathscr{B}\tensor
\mathscr{L}^{n}$ that does not vanish at $p\in X$.  
\end{proof}

\begin{rem}{Remark}{likeample}
Proposition~\ref{free} should be thought of as an analogue of the
statement that if $\mathscr{L}$ is an ample line bundle, then
$\mathscr{L}^n\tensor \mathscr{B}$ is globally generated for all $n\gg 0$.
{\nopagebreak \hfill $\Box$ \par \medskip}
\end{rem}

\begin{cor}\label{secgen}
For $k$ sufficiently large, the set theoretic base locus
of the linear system $|(2k-1)H-kE|$ on $\wtp$ is
$\P(\E)$. 
\end{cor}

\begin{proof}
Clearly, the base locus of $|(2k-1)H-kE|$ contains $\P(\E)$.
Note, however, that as $|2H-E|$ is base point free, the base locus of
$|(2k-1)H-kE|$ will stabilize for $k$ sufficiently large.  Hence it
suffices to show that if $p$ is a point not in $\P(\E)$,
then $|(2k-1)H-kE|$ is free at $p$ for all $k\gg 0$.  Now take
$\mathscr{B}=\O_{\wtp}(-H)$ and
$\mathscr{L}=\O_{\wtp}(2H-E)$ in Proposition~\ref{free}, and
use Theorem~\ref{embed}.  
\end{proof}

\begin{rem}{Notation}{notation}
For the rest of this section, write $\O(a,b,c)$ for
$\O_{\wt{M_2}}(aH+bE_1+cE_2)$, and write
$\mathscr{L}_k=\O_{\wt{M_2}}((2k-1)H-kE_1-E_2)$, $k\in \Q$.
\end{rem}

\begin{lemma}\label{nef}
$\mathscr{L}_k$ is nef for all $k$ sufficiently large. 
\end{lemma}

\begin{proof}
Letting $C\subset \wt{M_2}$ be an irreducible curve not contained
in $E_2$, we have $\mathscr{L}_k.C\geq 0$ for $k\gg 0$ by
Corollary~\ref{secgen}.  Letting $C'\subset E_2$, $\mathscr{L}_k\tensor 
\O_{E_2}$ is globally generated on $E_2$ for $k\gg 0$ by $(\ref{rest})$,
hence $\mathscr{L}_k.C'\geq 0$ and $\mathscr{L}_k$ is nef.  
\end{proof}

\begin{prop}\label{topofflip}
With hypotheses as above and for $k$ sufficiently large, 
the rational map on $\wt{M_2}$ induced by the linear
system $|\mathscr{L}_k|$ is a morphism, is an embedding off of $E_2$, and
its restriction to $E_2$ is the morphism $h:E_2\rightarrow \P(\F)$. 
\end{prop}

\begin{proof}
We first show that for $k$ sufficiently large $|\mathscr{L}_k|$
restricts to the complete linear system on $E_2$ associated to the
invertible sheaf $h^*f^*M^k\tensor \O_{E_2}(1)\tensor \pi^*\left(
\O_{\P(\E)}(-1)\right) $ (Cf. Lemma~\ref{itscomplete}).  For this it
suffices to prove $$H^1\left( \wt{M_2},\O(2k-1,-k,-2)\right)=0$$  

Writing $\mathscr{B}=\O(2k-1,-k,-2)$ and noting
$K_{\wt{M_2}}=\O(-n-1,n-r-1,n-2r-2)$:  
$$\mathscr{B}\tensor K^{-1}_{\wt{M_2}}\cong
\O(2k+n,-k-n+r+1,-n+2r)$$

Let $\alpha =\frac{k+n-r-1}{n-2r}$ and rewrite the right side as 
$\mathscr{L}_{\alpha}^{n-2r} \tensor \O(2,0,0)$.
For $k\gg 0$, $\mathscr{L}_{\alpha}$ is a nef $\Q$-divisor by
Lemma~\ref{nef}, hence $\mathscr{B}\tensor K^{-1}_{\wt{M_2}}$ is
big and nef and the vanishing holds by the Kawamata-Viehweg vanishing
theorem (note that by \cite[1.9]{mori}, a nef line bundle tensored with a
big and nef line bundle is again big and nef). 

To see $|\mathscr{L}_k|$ is a morphism, note that by
Corollary~\ref{secgen}, the support of the base scheme is contained in
$E_2$.  By what has just been proven, however, $|\mathscr{L}_k|$ has no
base points since the complete linear system on $E_2$ associated to
$h^*f^*M^k\tensor \O_{E_2}(1)\tensor \pi^*\left( \O_{\P(\E)}(-1)\right) $
induces a morphism.  Further, as $\O(2,-1,0)$ induces an embedding off
of $E_2$, $|\mathscr{L}_k|$ does as well for $k\gg 0$.
\end{proof}

\begin{rem}{Remark}{notsobadreally}
It is unfortunate that this proof gives no bound on $k$; however, there
is an important case when the value of $k$ can be determined.
Specifically, if $\SecX\subset \P^n$ is scheme theoretically defined by
cubics, then the line bundle $\mathscr{L}_2$ will be base point free and
$k=3$ suffices for Proposition~\ref{topofflip}.
{\nopagebreak \hfill $\Box$ \par \medskip}
\end{rem}

We simply write the morphism from Proposition~\ref{topofflip}
as $h:\wt{M_2}\rightarrow\P(|\mathscr{L}_k|)$.  Denote by $M_2$ the 
image variety of $h$.  Then $M_2\setminus
\P(\F)\cong \wtp\setminus \P(\E)$ by
Proposition~\ref{topofflip}.  

\begin{prop}\label{m2smooth}
$M_2$ is smooth.
\end{prop}

\begin{proof}
Note first that $M_2$ is normal since it is the image under a morphism
of a normal variety with reduced, connected fibers.

Let $Z=h^{-1}(p)\cong \P^1$, where $p\in
\P_{\H^2X}(\F)$. We have the normal bundle sequence:
$$0\rightarrow N_{Z/E_2}\rightarrow N_{Z/\wt{M_2}}\rightarrow
\O_{E_2}(E_2)\tensor \O_{Z}\rightarrow 0$$
Because $Z$ is a fiber of the $\P^1\times \P^{n-2r-1}$-bundle $E_2$ over
$\H^2X$ by Proposition~\ref{topofflip}, this sequence becomes:
$$\displaystyle 0\rightarrow \bigoplus_{1}^{n-2}\O_{\P^1}\rightarrow
N_{Z/\wt{M_2}}\rightarrow \O_{\P^1}(-1)\rightarrow 0$$

This sequence clearly splits and we apply a natural extension
of the smoothness portion of Castelnuovo's contractibility criterion
for surfaces (Cf. \cite[2.4]{aw}).  
\end{proof}

Because
$\wtp$ and $M_2$ are smooth and isomorphic in codimension one, we have
\cite[II.6.5]{hart} $\picard\wtp\cong \picard M_2$
and $H^0\left( \wtp, \O_{\wtp}(2H-E)\right) \cong
H^0\left( M_2, \O_{M_2}(2H-E)\right) $.
Therefore, the line bundle $\O_{M_2}(2H-E)$ induces a morphism
$f:M_2\rightarrow \P^s$ which is an embedding off of $\P(\F)$. 
Furthermore, because $M_2$ and $\P(\F)$ are smooth, \cite[1.1]{es}
implies that $h:\wt{M_2}\rightarrow M_2$ is the blow up
of $M_2$ along $\P(\F)$.

Collecting these results:

\begin{thm}\label{flip}
Let $(X,V)$ satisfy $(K_2)$ and assume $X\subset \P^n$ is smooth,
irreducible, and contains no lines or conics.  Then
there is a flip as pictured below with:
\begin{enumerate}
\item $\wtp$, $\wt{M_2}$, and $M_2$ smooth
\item $\wtp\setminus \P(\E)\cong M_2\setminus \P(\F)$
\item $h$ is the blow up of $M_2$ along $\P(\F)$
\item $\pi$ is the blow up of $\wtp$ along $\P(\E)$
\item $f$, induced by $\O_{M_2}(2H-E)$, is an embedding off
of $\P(\F)$, and the restriction of $f$ is
the projection $\P(\F)\rightarrow \H^2X$ 
\item $\wtv$, induced by $\O_{\wtp}(2H-E)$,
is an embedding off of $\P(\E)$, and the restriction of
$\wtv$ is the projection $\P(\E)\rightarrow \H^2X$
\end{enumerate}
\end{thm}

\begin{center}
{\begin{minipage}{1.5in}
\diagram
 & E_2 \dlto_{\pi} \drto^{h} & & & \wt{M_2} \dlto_{\pi}
\drto^{h} & \\
  \P(\E) \drto_{\wtv} & & \P(\F)
\dlto^f &  \wtp \drto_{\wtv} \dto & & M_2
\dlto^f \\
 & \H^2X & & \P^n \rdashed|>\tip_{\varphi }  & \P^s &
\enddiagram
\end{minipage}}
\end{center}
{\nopagebreak \hfill $\Box$ \par \medskip}

Note that the ample cone of the space $M_2$ is bounded on one side by
$\O_{M_2}(2H-E)$ and on the other by a line bundle of the form
$\O_{M_2}((2k-1)H-kE)$, $k\geq 2$, hence the assertion that this construction
yields a flip.  In fact, as $\O_{M_2}(3H-2E)$ is $f$-ample, this is a 
$(K+D)$-flip (in the sense of \cite[3.33]{km}) where
$D=(n+4)H-(n-r+1)E$.  When $r=1$ it is shown in \cite{bertram4}
that $D$ is log canonical, i.e. that this is a log flip.

\begin{rem}{Remark}{planequadscanbeok}
{\bf (on Conics)} The hypothesis that $X$ contain no conics 
is not always needed in order to construct this flip.  
The condition is imposed simply to keep the fibers of 
$\wtv:\wts\rightarrow\P^s$ equidimensional; problems occur 
when a variety has only a few conics.  In the case 
of quadratic Veronese embeddings where there is a unique plane 
quadric through any two points, it is not difficult to modify 
the results of the previous sections to achieve similar results.  
We do not do this, however, as this case is already understood from
the point of view of complete quadrics \cite{thad-coll},\cite{vain}.
{\nopagebreak \hfill $\Box$ \par \medskip}
\end{rem}

\begin{rem}{Example}{bigveroneses}
It is interesting to examine Theorem~\ref{flip} in cases where $\H^2X$
is well understood.  For example, in the case $X=\P^r$, it
is easy to see that $\H^2X$ is itself a $\P^2$-bundle over the Grassmannian
$\G(1,r)$ of lines in $\P^r$.  If we embed $\P^r$ via $\O(d)$, $d\geq
3$, $\wts$ has the particularly nice structure of a
$\P^1\times\P^2$-bundle over $\G(1,r)$ (in the case $d=2$, it is
simply a $\P^2$-bundle; the missing factor of $\P^1$ is due
precisely to the deficiency of $\SecX$).

Furthermore, as the ideal of $\SecX$ is generated by cubics
\cite{kanev}, the line bundle $\O_{M_2}(3H-2E)$ is globally generated
(see the discussion below).  Applying Kawamata-Viehweg vanishing
yields $$H^i(M_2, \O(kH-aE))=0, i>0, k>\frac{3}{2}(e+a-1)-(n+1)$$
where $n+1=$ ${r+d}\choose{r}$ (compare Corollary~\ref{littlevanishing}).
In the special case $k=2a-1$, this 
vanishing can be pulled directly back to $\wtp$.  More generally:
\end{rem}

\begin{cor}\label{secondvanishing}
With hypotheses as in Theorem~\ref{flip}, assume further that
$\SecX\subset\P^n$ is not a hypersurface and is cut out as a scheme by
cubics.  Then: $$H^i(\P^n,\I_X^a(2a-1))=0, i>0, a>n-3r-1$$ 
{\nopagebreak \hfill $\Box$ \par \medskip}
\end{cor}

The continuation of this process following \cite{thad} is
taken up in \cite{verm-flip2}.  We need to construct 
a birational morphism
$\wtv_2:M_2\rightarrow \P^{s_2}$ which contracts the image of $3$-secant
$2$-planes to points, and is an embedding off their union.  The natural
candidate for this morphism is the 
linear system associated to $\O_{M_2}(3H-2E)$, where we identify
$\picard\wtp\cong \picard M_2$.  Noting the fact
that $h^*\O_{M_2}(3H-2E)=\O_{\wt{M_2}}(3H-2E_1-E_2)$, it is not 
difficult to see (using Zariski's Main Theorem) that
this system will be globally generated if $\SecX\subset \P^n$ is
scheme theoretically defined by cubics.  

There are not yet theorems analogous to those for the quadric
generation of varieties.  However, there is evidence that such statements
should exist (Cf. \cite{kanev} where it is shown that
$Sec(v_d(\P^n))$ is ideal theoretically defined by cubics for all $d,n$).
Furthermore, the author proves set theoretic
statements for arbitrary smooth varieties in \cite{verm-flip2}.
These statements also contain information about the syzygies
among the generators that makes it possible to study the map
$\wtv_2$ in much the same way as $\wtv$ was studied above.

\end{article}

\end{document}